\documentclass{amsart}
\usepackage{amsfonts, amssymb}

\newtheorem{theorem}{Theorem}
\newtheorem{lemma}{Lemma}
\newtheorem{proposition}{Proposition}
\newtheorem{corollary}{Corollary}

\theoremstyle{definition}
\newtheorem{definition}{Definition}

\theoremstyle{remark}
\newtheorem{remark}{Remark}

\numberwithin{equation}{section}
\numberwithin{theorem}{section}
\numberwithin{lemma}{section}
\numberwithin{proposition}{section}
\numberwithin{corollary}{section}
\numberwithin{conjecture}{section}
\numberwithin{definition}{section}
\numberwithin{example}{section}
\numberwithin{xca}{section}
\numberwithin{remark}{section}
\numberwithin{claim}{section}
\numberwithin{question}{section}
\newcommand\lbb[1]{\label{#1}}
\newcommand{\thref}[1]{Theorem \ref{#1}}
\newcommand{\prref}[1]{Proposition \ref{#1}}
\newcommand{\leref}[1]{Lemma \ref{#1}}
\newcommand{\coref}[1]{Corollary \ref{#1}}

\newcommand{\deref}[1]{Definition \ref{#1}}

\newcommand{\reref}[1]{Remark \ref{#1}}

\newcommand{\seref}[1]{Section \ref{#1}}

\def\tt{\otimes}                               

\def\d{\partial}

\def\tint{{\textstyle\int}}

\def\res{|}                                    
\def\st{\; | \;}                               
\newcommand\nop[1]{{\rm{:}}#1{\rm{:}}}         
\def\vac{|0\rangle}                            

\def\ev{{\mathrm{ev}}}

\def\i{{\mathrm{i}}}     
\def\Cset{\mathbb{C}}       
\def\Zset{\mathbb{Z}}       
       
\def\al{\alpha}                         
\def\be{\beta}
\def\ga{\gamma}
\def\Ga{\Gamma}
\def\de{\delta}
\def\De{\Delta}
\def\ep{\varepsilon}
\def\eps{\epsilon}
\def\la{\lambda}
\def\La{\Lambda}

\def\Om{\Omega}

\def\si{\sigma}

\def\ze{\zeta}
\def\h{{\mathfrak{h}}}      

\def\A{\mathcal{A}}
\def\G{G}
\DeclareMathOperator{\rank}{rank}

\DeclareMathOperator{\Res}{Res}
\DeclareMathOperator{\Ind}{Ind}
\DeclareMathOperator{\ad}{ad}

\DeclareMathOperator{\id}{id}

\DeclareMathOperator{\Aut}{Aut}

\DeclareMathOperator{\End}{End}

\begin{document}

\title{Twisted Modules over Lattice Vertex Algebras}


\author{Bojko Bakalov
}

\address{Department of Mathematics, 
NCSU, Raleigh, NC 27695, USA}
\email{bojko\_bakalov@ncsu.edu}

\author{Victor G.~Kac
}

\address{Department of Mathematics, MIT, Cambridge, MA 02139, USA}
\email{kac@math.mit.edu}

\date{February 19, 2004}

\maketitle


\dedicatory{\textit{
Dedicated to Professor Ivan Todorov on the occasion of his 70th birthday.
}}

\section{Introduction}\lbb{sintro}
A vertex algebra is essentially the same as a chiral algebra
in conformal field theory \cite{BPZ, G}.
Vertex algebras arose naturally in the representation theory
of infinite-dimensional Lie algebras and in the construction
of the ``moonshine module'' for the Monster finite simple group
\cite{B1, FLM}.

If $V$ is a vertex algebra and $\Ga$ is a finite group 
of automorphisms of $V$, the subalgebra $V^\Ga$ of $\Ga$-invariant 
elements in $V$ is called an \emph{orbifold vertex algebra}.
Orbifolds play an important role in string theory;
in the physics literature they were introduced in
one of the earliest papers on conformal field theory 
\cite{DVVV, DHVW}. 
Recently, there have been numerous mathematical
papers on orbifolds. All these papers are concerned in some way or 
another with the problem of describing the
representations of $V^\Ga$ in terms of the
vertex algebra $V$ and the group $\Ga$.
However, the solution is known only in very special cases
and it is highly nontrivial.

Let $Q$ be an integral lattice. Then one can construct
a vertex algebra $V_Q$ called a \emph{lattice vertex algebra}
\cite{B1, FLM, K}. If $\bar\Ga$ is a finite group of isometries of $Q$,
its elements can be lifted to automorphisms of $V_Q$. One obtains
a group $\Ga\subset\Aut V_Q$, which is a central extension of $\bar\Ga$.
In \cite{BKT}, we construct a collection of $V_Q^\Ga$-modules,
and we compute their characters and modular transformations.
Some of these modules play an important role in the attempts of
a conformal field theory understanding of 
the fractional quantum Hall effect (see \cite{CGT}
and the references therein).

The present paper is the first step in the construction
of the $V_Q^\Ga$-modules \cite{BKT}.
Here we classify the so-called \emph{twisted $V_Q$-modules};
our main result is \thref{ttwvqm}.
In the case when $Q$ is a root lattice 
of a simple finite-dimensional Lie algebra 
and $\si$ is an element of its Weyl group,
our results agree with those of \cite{KP, KT}.
Some of our results were obtained independently in \cite{R},
and some special cases were studied earlier in \cite{D, X}.

{\flushleft{\bf{Acknowledgments.}}}
We are grateful to Ivan Todorov for many stimulating discussions
and for collaboration on \cite{BKT}.
We thank the organizers of the Varna Workshop
for inviting us to present our results
and for the inspiring workshop.
We also acknowledge the hospitality of the 
Erwin Schr\"odinger Institute,
where some of this work was done. 
The first author was supported in part by 
the Miller Institute for Basic Research in Science.
The second author was supported in part by 
NSF grant DMS-9970007.

\section{Preliminaries on Vertex Algebras}\lbb{sprva}

In this section we recall the definition of a 
vertex algebra and some of its properties, following the book~\cite{K}. 
Below $z,w,\dots$ denote formal commuting variables.
All vector spaces are over the field $\Cset$ of complex numbers.
\subsection{Local Fields}\lbb{slocf}
Let $V$ be a vector space.
A {\em field\/}  on $V$
is a formal power series in $z,z^{-1}$ of the form
\begin{equation}\lbb{field}
\phi(z)=\sum_{m\in\Zset}\phi_{(m)} z^{-m-1}, \qquad 
\phi_{(m)}\in\End(V),
\end{equation}
such that
\begin{equation}\lbb{gg}
\phi_{(m)} a = 0 \quad\text{for all }a\in V\text{ and }m\gg0.
\end{equation}
Note that all fields form a vector space invariant
under $\d_z$.

Two fields $\phi(z)$ and $\psi(z)$ are called {\em local\/}
with respect to each other if
\begin{equation}\lbb{lcl}
(z-w)^N\bigl[\phi(z),\psi(w)\bigr]=0 \quad\text{for }N\gg0
\end{equation}
as a formal power series in $z,z^{-1},w,w^{-1}$.
We will also say that the pair $(\phi,\psi)$ is local.
Obviously, if $(\phi,\psi)$ is local, then so are $(\psi,\phi)$
and $(\d\phi,\psi)$.

The locality \eqref{lcl} is equivalent to the commutator formula
\begin{equation}\lbb{lcl2}
\bigl[\phi(z),\psi(w)\bigr] =
\sum_{n\ge0} \bigl(\phi(w)_{(n)}\psi(w)\bigr) \d_w^{(n)}\de(z-w),
\end{equation}
where $\d_w^{(n)} = \d_w^n/n!$,
\begin{equation}\lbb{del}
\de(z-w)=\sum_{k\in\Zset} z^k w^{-k-1}
\end{equation}
is the formal $\de$-function, and
\begin{equation}\lbb{nprod}
\phi(w)_{(n)}\psi(w) = \Res_z (z-w)^n [\phi(z),\psi(w)\bigr], 
\qquad n\ge0,
\end{equation}
is the $n$th {\em product\/} of the fields $\phi$, $\psi$.
Here, as usual, $\Res_z$ denotes the coefficient at $z^{-1}$.
Note that the sum in \eqref{lcl2} is finite.

One can define $n$th product of fields for any $n\in\Zset$:
\begin{equation}\lbb{nprod2}
\begin{split}
\phi(w)_{(n)}\psi(w) &= 
 \Res_z \phi(z)\psi(w) \, i_{z,w} (z-w)^n
\\
&-\Res_z \psi(w)\phi(z) \, i_{w,z} (z-w)^n.
\end{split}
\end{equation}
Here $i_{z,w}$ means that we expand in the domain $|z|>|w|$, i.e.,
\begin{equation}\lbb{izw}
i_{z,w} (z-w)^n = \sum_{k=0}^\infty \binom{n}{k} z^{n-k} (-w)^k,
\end{equation}
while
\begin{equation}\lbb{iwz}
i_{w,z} (z-w)^n = \sum_{k=0}^\infty \binom{n}{k} z^k (-w)^{n-k}.
\end{equation}
In particular, $\de(z-w) = i_{z,w}(z-w)^{-1} - i_{w,z}(z-w)^{-1}$.

The $(-1)$st product is called the {\em normally ordered product\/}
and is denoted by $\nop{\phi\psi}$. One has
\begin{equation}\lbb{nop1}
\nop{\phi(w)\psi(w)} = \phi(w)^+\psi(w) + \psi(w)\phi(w)^-,
\end{equation}
where
\begin{equation}\lbb{nop2}
\phi(w)^+ = \sum_{m<0}\phi_{(m)} w^{-m-1}, \quad 
\phi(w)^- = \sum_{m\ge0}\phi_{(m)} w^{-m-1}.
\end{equation}
One also has
\begin{equation}\lbb{nprod3}
\phi(w)_{(-n-1)}\psi(w) = \nop{(\d_w^{(n)}\phi(w))\psi(w)},
\qquad n\ge0.
\end{equation}


\subsection{Vertex Algebras}\lbb{sdva}
\begin{definition}\lbb{dva}
A {\em vertex algebra\/} is a vector space $V$
(space of states) endowed with a vector
$\vac\in V$ (vacuum vector), an endomorphism
$T$ (infinitesimal translation operator)
and a linear map {}from $V$ to the space of fields on $V$
(state-field correspondence)
\begin{equation}\lbb{y}
a\mapsto Y(a,z)=\sum_{m\in\Zset}a_{(m)}z^{-m-1},
\end{equation}
such that the following properties hold:

locality axiom: all fields $Y(a,z)$ are local with respect to each other,

translation covariance:
\begin{equation}\lbb{t}
\bigl[T,Y(a,z)\bigr]=\d_z Y(a,z),
\end{equation}

vacuum axioms:
\begin{align}
\lbb{vac1}
&Y(\vac,z)=\id_V, \;\;\; T\vac=0,
\\
\lbb{vac2}
&Y(a,z)\vac-a \in zV[[z]].
\end{align}
\end{definition}

Here are some corollaries of the definition:
\begin{equation}\lbb{taa}
Ta=a_{(-2)}\vac, \quad
Y(Ta,z)=\d_zY(a,z).
\end{equation}
We also have the skew-symmetry
\begin{equation}\lbb{sks}
Y(a,z)b=e^{zT}Y(b,-z)a.
\end{equation}
The most important property of a vertex algebra $V$ is the following
{\em Borcherds identity\/} (which along with the vacuum axioms
provides an equivalent definition of a vertex algebra):
\begin{equation}\lbb{bor}
\begin{split}
\Res_{z-w} Y(Y(a,z-w)b,w) \, i_{w,z-w}&F(z,w)
\\
=\Res_z Y(a,z)Y(b,w) \, i_{z,w}&F(z,w)
\\
-\Res_z Y(b,w)Y(a,z) \, i_{w,z}&F(z,w)
\end{split}
\end{equation}
for any $a,b\in V$ and any rational function $F(z,w)$ with poles only at 
$z=0$, $w=0$ or $z=w$.

Let us give some consequences of the Borcherds identity.
Taking $F(z,w)=(z-w)^n\de(u-z)$, viewed as a series in 
$u,u^{-1}$, we obtain
\begin{equation}\lbb{bor2}
\begin{split}
Y(a,u)&Y(b,w) \, i_{u,w}(u-w)^n 
-Y(b,w)Y(a,u) \, i_{w,u}(u-w)^n
\\
&=\sum_{m=0}^\infty Y(a_{(m+n)}b,w) \, \d_w^{(m)}\de(u-w).
\end{split}
\end{equation}
For $n=0$ this gives the commutator formula
\begin{equation}\lbb{comm}
\bigl[Y(a,u),Y(b,w)\bigr] =
\sum_{m=0}^\infty Y(a_{(m)}b,w) \, \d_w^{(m)}\de(u-w),
\end{equation}
(which implies locality). Taking $\Res_u$ of \eqref{bor2}, we get
\begin{equation}\lbb{npr}
Y(a,w)_{(n)}Y(b,w)=Y(a_{(n)}b,w), \qquad n\in\Zset.
\end{equation}


\subsection{Conformal Vertex Algebras}\lbb{sconfdva}
\begin{definition}\lbb{dconfva}
A vertex algebra $V$ is {\em conformal\/} of central charge (or rank)
$c\in\Cset$ if there exists a vector $\nu\in V$
(the {conformal vector}) with the following properties:

(i) The modes $L_n$, $n\in\Zset$, of the field
$L(z) \equiv Y(\nu,z)=\sum_{n\in\Zset}L_n z^{-n-2}$
give a representation of the Virasoro algebra with
central charge $c$. (The field $L(z)$ is called
the energy-momentum field in CFT.)

(ii) $L_{-1}$ is the infinitesimal translation operator $T$.

(iii) The operator $L_0$ is diagonalizable with a non-negative spectrum.
($L_0$ is called the energy operator or Hamiltonian.)
\end{definition}
%

In a conformal vertex algebra $V$,  \eqref{comm} implies
\begin{equation}\lbb{l02}
[L_0,Y(a,z)] = zY(L_{-1}a,z) + Y(L_0a,z).
\end{equation}
The operator $L_0$ defines a gradation of $V$,
$V=\bigoplus_{\De\ge0}V(\De)$, such that
\begin{equation}\lbb{l0}
L_0\res_{V(\De)} = \De\id_{V(\De)}.
\end{equation}
For $a\in V(\Delta)$, the field $Y(a,z)$ is of {\em conformal weight\/}
$\Delta$, i.e., $\deg a_{(m)}=\Delta-m-1$ for all $m\in\Zset$.
The field $L(z)$ is of conformal weight $2$ and $\nu=L_{-2}\vac$.

\section{Twisted Modules over Vertex Algebras}\lbb{stwmod}
In this section we study twisted modules over vertex algebras 
(cf.\ \cite{DL, FLM}).
It seems that some of our results are new even
in the untwisted case.

\subsection{Definition in Terms of Borcherds Identity}\lbb{sdbid}
Let $V$ be a vertex algebra and $\si$ be an automorphism of $V$
of finite order $N$. We let $\eps=e^{2\pi\i/N}$ and
$V_{j} = \{a\in V \st \si a=\eps^{-j}a \}$, $0\le j\le N-1$.

An {\em $N$-twisted field\/} $\phi(z)$
on a vector space $M$ is a formal power series in $z^{1/N},z^{-1/N}$ 
of the form
\begin{equation}\lbb{twfield}
\phi(z)=\sum_{m\in\frac1N\Zset}\phi_{(m)} z^{-m-1}, \qquad 
\phi_{(m)}\in\End(M),
\end{equation}
such that
\begin{equation}\lbb{twgg}
\phi_{(m)} v = 0 \quad\text{for all }v\in M\text{ and }m\gg0.
\end{equation}
For any integer $k$, we will denote by $\phi(e^{2\pi\i k} z)$ the field
obtained from $\phi(z)$ by substituting $z^{1/N}$ with $\eps^k z^{1/N}$, i.e.,
$\phi(e^{2\pi\i k} z) 
= \sum_{m\in\frac1N\Zset}\phi_{(m)} e^{-2\pi\i km} z^{-m-1}$.

\begin{definition}\lbb{dvrep}
A {\em $\si$-twisted $V$-module\/} is a
vector space $M$ endowed with 
a linear map {}from $V$ to the space of $N$-twisted fields on $M$,
\begin{align}
\lbb{ym}
a\mapsto Y^M(a,z)&=\sum_{m\in\frac1N\Zset}a^M_{(m)}z^{-m-1}, 
\\
\intertext{such that for all $a\in V$}
\lbb{ginv}
Y^M(\si a,z) &= Y^M(a,e^{2\pi\i}z),
\\
\lbb{vacm}
Y^M(\vac,z) &= \id_M, 
\end{align}
and the following twisted Borcherds identity holds for any 
$a\in V_j$, $b\in V$, $0\le j\le N-1$,
and any rational function $F(z,w)$ with poles 
only at $z=0$, $w=0$ or $z=w$:
\begin{equation}\lbb{twbor}
\begin{split}
\Res_{z-w} Y^M(Y(a,z-w)b,w) \, i_{w,z-w} \, z^{j/N} &F(z,w)
\\
=\Res_z Y^M(a,z)Y^M(b,w) \, i_{z,w} \, z^{j/N} &F(z,w)
\\
-\Res_z Y^M(b,w)Y^M(a,z) \, i_{w,z} \, z^{j/N} &F(z,w).
\end{split}
\end{equation}

Of course, for $\si=1$, a $1$-twisted module is called just a {\em module\/}.
\end{definition}
\begin{remark}\lbb{rtwbor}
The Borcherds identity \eqref{twbor} is equivalent to the following
collection of identities ($a\in V_j$, $b\in V$, $c\in M$,
$m\in\frac jN+\Zset$, $n\in\Zset$, $k\in\frac1N\Zset$):
\begin{equation}\lbb{twboco}
\begin{split}
\sum_{i=0}^\infty \binom{m}{i} (a_{(n+i)}b)^M_{(m+k-i)}c
= \sum_{i=0}^\infty (-1)^i &\binom{n}{i} a^M_{(m+n-i)}(b^M_{(k+i)}c)
\\
- \sum_{i=0}^\infty (-1)^{i+n} &\binom{n}{i} b^M_{(k+n-i)}(a^M_{(m+i)}c).
\end{split}
\end{equation}
\end{remark}

When $V$ is a conformal vertex algebra we will assume that the
automorphism $\si$ preserves the conformal vector: $\si(\nu)=\nu$.
Then, for any $\si$-twisted $V$-module $M$, the modes 
$L_n^M$ of
$Y^M(\nu,z)=\sum_{n\in\Zset}L_n^M z^{-n-2}$
give a representation in $M$ of the Virasoro algebra with
central charge $c=\rank V$.

\begin{definition}\lbb{dvrep2}
A {\em $\si$-twisted module over a conformal vertex algebra $V$\/} 
is a $\si$-twisted $V$-module $M$ in the sense of \deref{dvrep}
satisfying the additional
requirement that the operator $L_0^M$ 
be diagonalizable with finite-dimensional eigenspaces.
\end{definition}

\subsection{Consequences of the Definition}\lbb{subconsd}
Let $M$ be a $\si$-twisted $V$-module.
In this subsection we will derive some consequences of \deref{dvrep}.

First of all, note that, by \eqref{ginv}, $Y^M(a,z)z^{j/N}$
contains only integer powers of $z$ for $a\in V_j$.
For $0\le j\le N-1$, we introduce the $N$-twisted $\de$-functions
\begin{equation}\lbb{twdel}
\de_j(z-w) = z^{j/N}w^{-j/N}\de(z-w) 
= \sum_{k\in\frac jN+\Zset} z^k w^{-k-1} . 
\end{equation}
Then 
\begin{equation}\lbb{twdel2}
\Res_z Y^M(a,z) \de_j(w-z) = Y^M(\pi_j a,w),
\end{equation}
where 
\begin{equation}\lbb{pij}
\pi_j = \frac1N \sum_{k=0}^{N-1} \eps^{kj} \si^k
\end{equation}
is the projection onto $V_{j}$. 

Putting $F(z,w)=z^{-j/N}(z-w)^n\de_j(u-z)$ in the twisted Borcherds identity
\eqref{twbor}, we obtain a twisted version of \eqref{bor2}:
\begin{equation}\lbb{twbor2}
\begin{split}
Y^M&(a,u)Y^M(b,w) \, i_{u,w}(u-w)^n 
-Y^M(b,w)Y^M(a,u) \, i_{w,u}(u-w)^n
\\
&=\sum_{m=0}^\infty Y^M(a_{(m+n)}b,w) \, \d_w^{(m)}\de_j(u-w),
\qquad a\in V_j, n\in\Zset.
\end{split}
\end{equation}
The collection of identities \eqref{twbor2} for all $n\in\Zset$
is equivalent to \eqref{twbor}.
Note that, as in the untwisted case, all the fields $Y^M(a,z)$
are local with respect to each other, since, letting $n=0$ in \eqref{twbor2},
we get the twisted commutator formula $(a\in V_j)$:
\begin{equation}\lbb{twcomm}
\bigl[Y^M(a,u),Y^M(b,w)\bigr] =
\sum_{m=0}^\infty Y^M(a_{(m)}b,w) \, \d_w^{(m)}\de_j(u-w) .
\end{equation}

Let us define the $n$th product ($n\in\Zset$) of two twisted
fields $Y^M(a,w)$, $Y^M(b,w)$ for $a\in V_j$ as 
$\Res_u u^{j/N}w^{-j/N}$ of the left hand side of \eqref{twbor2}.
For $n=-1$, this definition coincides with 
\eqref{nop1}, \eqref{nop2}, but for $n<-1$ it differs from \eqref{nprod3}.
By \eqref{twbor2}, we have for $a\in V_j$
\begin{equation}\lbb{twnpr}
Y^M(a,w)_{(n)}Y^M(b,w) = 
\sum_{m=0}^\infty \binom{j/N}{m} w^{-m} \, Y^M(a_{(m+n)}b,w) . 
\end{equation}
In particular, letting $b=\vac$ and $n=-2$ in \eqref{twnpr}, we get
\begin{equation}\lbb{twtaa}
Y^M(Ta,z)=\d_zY^M(a,z).
\end{equation}

When $V$ is a conformal vertex algebra,
the commutator formula \eqref{twcomm} implies:
\begin{align}
\lbb{tm}
[T^M,Y^M(a,z)] &= \d_z Y^M(a,z), \\
\lbb{l0m}
[L_0^M,Y^M(a,z)] &= z\d_z Y^M(a,z) + Y^M(L_0^M a,z),
\end{align}
where $T^M := L_{-1}^M$.

\subsection{Associativity of Twisted Fields }\lbb{sass}
We will try to invert formula \eqref{twnpr}, i.e., to
express $Y^M(a_{(n)}b,w)$ in terms of $Y^M(a,w)_{(n)}Y^M(b,w)$.
To this end, we multiply both sides of \eqref{twnpr} by 
$z^{-n-1}w^{j/N}$ and sum over $n\in\Zset$. Using the properties
of the $\de$-function and Taylor's formula, we get
\begin{align}
\notag
i_{w,z} &(z+w)^{j/N} \, Y^M(Y(a,z)b,w) 
= i_{z,w} (z+w)^{j/N} \, Y^M(a,z+w)Y^M(b,w)
\\
\lbb{as1}
&- \sum_{k=0}^\infty Y^M(b,w) a^M_{(k+j/N)} \,  (-\d_w)^{(k)}\de(z-(-w)),
\qquad a\in V_j.
\end{align}
Of course, we cannot divide \eqref{as1} by $(z+w)^{j/N}$. Note, however, that
the sum in \eqref{as1} becomes finite when applied to any $v\in M$.
Using that $(z+w)^n (-\d_w)^{(k)}\de(z-(-w)) = 0$ for $n>k$,
we obtain the {\em associativity\/} of twisted fields
\begin{equation}\lbb{as2}
\begin{split}
i_{w,z} (z&+w)^{n+j/N} \, Y^M(Y(a,z)b,w)v
\\
&= i_{z,w} (z+w)^{n+j/N} \, Y^M(a,z+w)Y^M(b,w)v,
\\
&\qquad\qquad\qquad a\in V_j, b\in V, v\in M, n\gg0.
\end{split}
\end{equation}
Here we can take the minimal $n\ge0$
such that $a^M_{(k+j/N)}v=0$ for $k\ge n$.

Note that the powers of $z$ and $w$ in both sides of \eqref{as2}
are bounded from below. Therefore, we can multiply both sides
by $i_{w,z}(z+w)^{-n-j/N}$ to get
($a\in V_j, b\in V, v\in M$):
\begin{equation}\lbb{as3}
\begin{split}
&Y^M(Y(a,z)b,w)v
\\
&= (i_{w,z}(z+w)^{-n-j/N}) \, i_{z,w} (z+w)^{n+j/N} \, Y^M(a,z+w)Y^M(b,w)v.
\end{split}
\end{equation}
%

\subsection{Definition in Terms of Associativity}\lbb{sdass}
We will show that the associativity \eqref{as2}, together with
\eqref{ginv}, \eqref{vacm}, 
implies the twisted Borcherds identity \eqref{twbor}.

First let $b=\vac$ in \eqref{as2}. Using \eqref{vacm} and the identity
$Y(a,z)\vac = e^{zT}a$ (cf.\ \eqref{sks}), we get:
\begin{equation*}
i_{w,z} (z+w)^{n+j/N} \, Y^M(e^{zT}a,w)v
= i_{z,w} (z+w)^{n+j/N} \, Y^M(a,z+w)v.
\end{equation*}
Notice that the right-hand side contains only non-negative integer powers
of $z+w$.
This implies $Y^M(e^{zT}a,w) = i_{w,z} Y^M(a,z+w)$ and, in particular,
Eq.~\eqref{twtaa}.

Now let $a\in V_j, b\in V_k, v\in M$, and let $n,p\in\Zset$ be such that
$a^M_{(m+j/N)}v=0$ for $m\ge n$, $b^M_{(m+k/N)}v=0$ for $m\ge p$.
Replacing in \eqref{as2} $b$ with $e^{uT}b$, we obtain:
\begin{equation}\lbb{as4}
\begin{split}
i_{z,w} i_{w,u} & (z+w)^{n+j/N} \, Y^M(a,z+w) Y^M(b,u+w) v
\\
&= i_{w,z} i_{w,u} i_{z,u} (z+w)^{n+j/N} \, Y^M(Y(a,z-u)b,u+w) v.
\end{split}
\end{equation}
Note that, if we multiply the left-hand side by $i_{w,u} (u+w)^{p+k/N}$,
it will contain only non-negative integer powers of $u+w$ and hence of $w$.
Therefore it makes sense to put $w=0$:
\begin{equation}\lbb{as5}
\begin{split}
z^{n+j/N} u^{p+k/N} \, & Y^M(a,z) Y^M(b,u) v
\\
= & i_{w,z} i_{w,u} i_{z,u} (z+w)^{n+j/N} (u+w)^{p+k/N} 
\\
&\times
Y^M(Y(a,z-u)b,u+w) v \big|_{w=0}.
\end{split}
\end{equation}
Interchanging the roles of $a$ and $b$ and using \eqref{sks}, we get:
\begin{equation}\lbb{as6}
\begin{split}
z^{n+j/N} u^{p+k/N} \, & Y^M(b,u) Y^M(a,z) v
\\
= & i_{w,z} i_{w,u} i_{u,z} (z+w)^{n+j/N} (u+w)^{p+k/N}
\\
&\times
Y^M(Y(a,z-u)b,u+w) v \big|_{w=0}.
\end{split}
\end{equation}

Notice that 
\begin{equation*}
i_{w,z} (z+w)^{n+j/N} = i_{w,u} \sum_{i\in\Zset_+} \binom{n+j/N}{i}
(u+w)^{n-i+j/N} (z-u)^i
\end{equation*}
and 
\begin{equation*}
i_{z,u} (z-u)^{-m-1} - i_{u,z} (z-u)^{-m-1} = \d_u^{(m)} \de(z-u),
\qquad m\ge0.
\end{equation*}
Using this and \eqref{as5}, \eqref{as6}, we get for $l\in\Zset$:
\begin{equation*}
\begin{split}
i_{z,u} (z-u)^l & z^{n+j/N} u^{p+k/N} \, Y^M(a,z) Y^M(b,u) v
\\
&- i_{u,z} (z-u)^l z^{n+j/N} u^{p+k/N} \, Y^M(b,u) Y^M(a,z) v
\\
=i_{w,u} & \sum_{\substack{ i,m\in\Zset \\ i\ge0, m\ge i+l }}
\binom{n+j/N}{i} (u+w)^{n+p-i+(j+k)/N}
\\
&\times
Y^M(a_{(m)}b,u+w) v \, \d_u^{(m-i-l)} \de(z-u) \big|_{w=0}.
\end{split}
\end{equation*}
{}From here it is easy to deduce \eqref{twbor2}, which implies \eqref{twbor}.

Therefore, we have the following equivalent definition of a 
$\si$-twisted $V$-module.
\begin{proposition}\lbb{pdefas}
A $\si$-twisted $V$-module is the same as a vector space $M$ endowed with 
a linear map \eqref{ym}
{}from $V$ to the space of $N$-twisted fields on $M$,
satisfying \eqref{ginv}, \eqref{vacm} and \eqref{as2}.
\end{proposition}

\section{Twisted Modules over a Lattice Vertex Algebra}\lbb{stmvq}
In the first subsection we introduce the main object of our study:
the lattice vertex algebra $V_Q$. The remainder of the section
is devoted to the classification of all $\si$-twisted $V_Q$-modules,
where $\si$ is an automorphism of the lattice $Q$.

\subsection{Lattice Vertex Algebras}\lbb{slva}
The purpose of this subsection is to fix the notation
and review some properties of lattice vertex algebras.
For more details, see~\cite{K}.

Let $Q$ be an integral lattice of rank $l$.
We denote the bilinear form on $Q$ by $(\cdot|\cdot)$, and write
$|\al|^2=(\al|\al)$ for $\al\in Q$. We extend the bilinear form
to $\h = \Cset\tt_\Zset Q$ by $\Cset$-bilinearity.
There exists a bimultiplicative function $\ep\colon Q\times Q\to\{\pm1\}$
satisfying
\begin{equation}\lbb{epal}
\ep(\al,\al) = (-1)^{|\al|^2(|\al|^2+1)/2}, \qquad \al\in Q.
\end{equation}
Then by bimultiplicativity
\begin{equation}\lbb{epalbe}
\ep(\al,\be)\ep(\be,\al) = (-1)^{(\al|\be) + |\al|^2|\be|^2},
\qquad \al,\be\in Q.
\end{equation}
Introduce the twisted group algebra $\Cset_\ep[Q]$: 
it has a basis $\{e^\al\}_{\al\in Q}$ over $\Cset$ and multiplication
\begin{equation}\lbb{ealebe}
e^\al e^\be = \ep(\al,\be) e^{\al+\be},
\qquad \al,\be\in Q.
\end{equation}

Let $\hat\h = \h[t,t^{-1}]\oplus\Cset K$ be the 
{\em Heisenberg current algebra};
this is a Lie algebra with the bracket
\begin{equation}\lbb{htht}
[ht^m, h't^n] = m\de_{m,-n}(h|h')K, \quad [ht^m,K]=0,
\qquad h,h'\in\h.
\end{equation}
It has a unique irreducible representation of level $1$ (i.e., with $K=1$)
on the Fock space $S=S(\h[t^{-1}]t^{-1})$ such that $\h[t^{-1}]t^{-1}$
acts by multiplication and $\h[t]1=0$. This representation
extends to the space $V_Q=S\tt\Cset_\ep[Q]$ by
\begin{equation}\lbb{hteal}
(ht^m)(s\tt e^\al) = \bigl( ht^m + \de_{m,0}(h|\al) \bigr) s\tt e^\al 
\qquad\text{for}\;\; m\ge0.
\end{equation}

We define a representation of the algebra $\Cset_\ep[Q]$ in $V_Q$
by left multiplication:
\begin{equation}\lbb{egaeal1}
e^\ga(s\tt e^\al) = \ep(\ga,\al)s\tt e^{\al+\ga}.
\end{equation}
This gives rise to a representation in $V_Q$ of the associative
algebra $\A_Q = U(\hat\h)\tt\Cset_\ep[Q]$, which is a 
``twisted'' tensor product of the universal enveloping algebra $U(\hat\h)$
of $\hat\h$ and the algebra $\Cset_\ep[Q]$ by the relation
\begin{equation}\lbb{ealhtm}
e^\al(ht^m) = (ht^m-\de_{m,0}(\al|h)) e^\al.
\end{equation}
Here and further $e^\al$ (respectively $ht^m$) stands for
$1\tt e^\al$ (respectively $ht^m \tt1$).
The algebra $\A_Q$ has a $\Zset_2$-gradation (i.e., it is an associative
superalgebra) defined by
\begin{equation}\lbb{pueal1}
p(u\tt e^\al) = |\al|^2 \mod 2\Zset.
\end{equation}
This induces a $\Zset_2$-gradation on $V_Q$:
\begin{equation}\lbb{pueal2}
p(s\tt e^\al) = |\al|^2 \mod 2\Zset.
\end{equation}

Introduce the following fields on $V_Q$ (called {\em currents}):
\begin{equation}\lbb{hz}
h(z) = \sum_{m\in\Zset} (ht^m)z^{-m-1},
\qquad h\in\h.
\end{equation}
Then the commutation relations \eqref{htht} for $K=1$ can be rewritten as
\begin{equation}\lbb{hzhz}
[h(z),h'(w)] = (h|h')\,\d_w\de(z-w),
\qquad h,h'\in\h.
\end{equation}
Hence all the fields $h(z)$ are local with respect to each other.

For $\al\in Q$, introduce the {\em vertex operator\/}
\begin{equation}\lbb{gaal}
\begin{split}
Y_\al(z) &= e^\al \,\nop{\exp\tint\al(z)}
\\
&\equiv e^\al z^\al \exp\Bigl( \sum_{n<0} (\al t^n)\frac{z^{-n}}{-n} \Bigr)
               \exp\Bigl( \sum_{n>0} (\al t^n)\frac{z^{-n}}{-n} \Bigr)
\end{split}
\end{equation}
of parity $|\al|^2 \mod 2\Zset$.
The vertex operator $Y_\al(z)$ is a field on $V_Q$, and it
is local with respect to $h(z)$ because
\begin{equation}\lbb{gaalhz}
[h(z),Y_\al(w)] = (h|\al)Y_\al(w)\de(z-w),
\qquad h\in\h,\al\in Q.
\end{equation}
Moreover, the fields $Y_\al(z)$ are local among themselves.
This follows from \eqref{epalbe} and the formula
\begin{equation}\lbb{gaalgabe}
Y_\al(z)Y_\be(w) = \ep(\al,\be) (z-w)^{(\al|\be)} \, e^{\al+\be} \,
\nop{\exp\bigl( \tint\al(z)+\tint\be(w) \bigr)} \,.
\end{equation}

\begin{theorem}\lbb{pvq}
The fields $Y(ht^{-1},z)=h(z)$ $(h\in\h)$, of parity $0$,
and $Y(e^\al,z)=Y_\al(z)$ $(\al\in Q)$, of parity 
$p(e^\al)=|\al|^2 \mod 2\Zset$, generate a vertex algebra structure
on $V_Q=S\tt\Cset_\ep[Q]$ with the vacuum vector $\vac=1\tt1$ and
the operator $T$ defined by
\begin{equation}\lbb{tvq}
[T,ht^m]=-m \, ht^{m-1}, \quad Te^\al=(\al t^{-1})e^\al,
\qquad h\in\h,\al\in Q.
\end{equation}
This vertex algebra
is conformal of rank $l=\rank Q$ with the conformal vector
\begin{equation}\lbb{nuvq}
\nu = \frac12 \sum_{i=1}^l (a^it^{-1})(b^it^{-1}),
\end{equation}
where $\{a^i\}$, $\{b^i\}$ are dual bases of\/ $\h$.
\end{theorem}

Let $\si$ be an automorphism of the lattice $Q$ of finite order $N$.
For $0\le j\le N-1$, let
\begin{equation}\lbb{hj}
\h_j = \{h\in\h \st \si h=\eps^{-j}h \}, 
\qquad \eps=e^{2\pi\i/N}.
\end{equation}
Then $\h_j$ and $\h_k$ are orthogonal unless $j+k\equiv 0\mod N\Zset$.
Since there exists a unique, up to equivalence, $\pm1$-valued
$2$-cocycle $\ep$ on $Q$ satisfying \eqref{epal}, the cocycles
$\ep(\al,\be)$ and $\ep(\si\al,\si\be)$ are equivalent.
Hence there exists a function $\eta\colon Q\to\{\pm1\}$ such that
\begin{equation}\lbb{etaep}
\eta(\al)\eta(\be)\ep(\al,\be) = \eta(\al+\be)\ep(\si\al,\si\be),
\qquad \al,\be\in Q.
\end{equation}
Moreover, $\eta$ can be chosen in such a way that
\begin{equation}\lbb{etah0}
\eta(\al)=1 \quad\text{for}\;\; \al\in Q\cap\h_0.
\end{equation}
\begin{proposition}\lbb{pautvq}
Any automorphism $\si$ of $Q$ can be lifted to an automorphism of 
the vertex algebra $V_Q$ so that
\begin{equation}\lbb{autvq}
\si(ht^m)=\si(h)t^m, \quad \si(e^\al)=\eta(\al)^{-1}e^{\si\al},
\qquad h\in\h,\al\in Q.
\end{equation}
It fixes the conformal vector $\nu${\rm:} $\si(\nu)=\nu$.
\end{proposition}
\begin{proof}
Follows from the observation that $\si$ defines an automorphism
of the associative superalgebra $\A_Q$.
\end{proof}
\begin{remark}\lbb{rordsi}
If $\si\colon Q\to Q$ is an automorphism of order $N$,
then its lifting $\si\colon V_Q\to V_Q$ defined by \eqref{autvq}
has order $N$ or $2N$.
\end{remark}

\subsection{Twisted Vertex Operators}\lbb{sstwvo}
Let $\si$ be an automorphism of $Q$ lifted to an automorphism of $V_Q$,
as in \prref{pautvq}. 
We use the notation from \seref{stwmod}
for $V=V_Q$. 
Let $M$ be  a $\si$-twisted $V_Q$-module (see \deref{dvrep2}).
We will study the fields
$Y^M(ht^{-1},z)$ and $Y^M(e^\al,z)$ for $h\in\h,\al\in Q$.

Comparing the commutator formulas \eqref{comm} and \eqref{twcomm},
we see that Eqs.\ \eqref{hzhz}, \eqref{gaalhz} immediately imply
\begin{align}
\lbb{twc1}
[Y^M(ht^{-1},z),Y^M(h't^{-1},w)] &= (h|h')\,\d_w\de_j(z-w),
\\
\lbb{twc2}
[Y^M(ht^{-1},z),Y^M(e^\al,w)] &= (h|\al)\, Y^M(e^\al,w)\de_j(z-w),
\\
\notag
&\qquad\qquad h\in\h_j,h'\in\h, \; \al\in Q.
\end{align}
These formulas can be restated as:
\begin{align}
\lbb{twc3}
[h^M_{(m)},(h')^M_{(n)}] &= (\pi_{Nm} h|h') \, m\de_{m,-n},
\\
\lbb{twc4}
[h^M_{(m)},Y^M(e^\al,w)] &= (\pi_{Nm} h|\al) \, w^m Y^M(e^\al,w),
\\
\notag
&\qquad\qquad h,h'\in\h, \; m,n\in\tfrac1N\Zset, \; \al\in Q,
\end{align}
where $\pi_j$ is the projection of $\h$ onto $\h_{j\!\!\!\mod N\Zset}$.
(Note that this $\pi_j$ is the restriction to $(\h t^{-1})\vac$ 
of the $\pi_j$ defined by \eqref{pij}). 
In the sequel we will use the notation $h_0=\pi_0 h$ for $h\in\h$.

Also note that, for $a=ht^{-1}$, \eqref{ginv} is equivalent to
\begin{equation}\lbb{twc5}
(\si h)^M_{(m)} = h^M_{(m)} \, e^{-2\pi\i m},
\qquad h\in\h, \; m\in\tfrac1N\Zset.
\end{equation}

\begin{lemma}\lbb{lymea}
There exist operators $U^M_\al$ $(\al\in Q)$ on $M$ such that
\begin{equation}\lbb{ymeal}
Y^M(e^\al,z) =  z^{b_\al} U^M_\al E^M_\al(z)
\end{equation}
where
\begin{equation}\lbb{dal}
b_\al = (|\al_0|^2 - |\al|^2)/2
\end{equation}
and
\begin{equation}\lbb{emal}
\begin{split}
E^M_\al(z) &= \nop{\exp\tint Y^M(\al t^{-1},z)}
\\
&\equiv z^{ \al^M_{(0)} } \,
\exp\Bigl( \sum_{n\in\frac1N\Zset_{<0}}
           \al^M_{(n)}\frac{z^{-n}}{-n} \Bigr)
\exp\Bigl( \sum_{n\in\frac1N\Zset_{>0}}
           \al^M_{(n)}\frac{z^{-n}}{-n} \Bigr).
\end{split}
\end{equation}
The operators $U^M_\al$ satisfy
\begin{align}
\lbb{hual}
[h^M_{(m)},U^M_\al] &= \de_{m,0}(h_0|\al) \, U^M_\al,
\qquad h\in\h, \; m\in\tfrac1N\Zset,
\\
\lbb{gual}
U^M_{\si\al} &= \eta(\al) U^M_\al e^{ 2\pi\i(b_\al+\al^M_{(0)}) }.
\end{align}
\end{lemma}
\begin{proof}
Define the operators 
\begin{equation*}
U^M_\al(z) = \exp\Bigl( \sum_{n\in\frac1N\Zset_{<0}}
           \al^M_{(n)}\frac{z^{-n}}{n} \Bigr) 
\, Y^M(e^\al,z) \, z^{ -\al^M_{(0)} }
\exp\Bigl( \sum_{n\in\frac1N\Zset_{>0}}
           \al^M_{(n)}\frac{z^{-n}}{n} \Bigr).
\end{equation*}
Then, by \eqref{twc3}, \eqref{twc4},
\begin{equation*}
[h^M_{(m)},U^M_\al(z)] = \de_{m,0}(h_0|\al) \, U^M_\al(z),
\qquad h\in\h.
\end{equation*}
In particular, $U^M_\al(z)$ commutes with $\al^M_{(n)}$
for $n<0$; hence, we have
\begin{equation}\lbb{ymeal2}
Y^M(e^\al,z) =  U^M_\al(z) E^M_\al(z).
\end{equation}

We will deduce Eqs.\ \eqref{ymeal} and \eqref{hual} once we show
that $U^M_\al := z^{ -b_\al} U^M_\al(z)$
is independent of $z$.
To this end we use the translation invariance \eqref{twtaa}.
By \eqref{tvq}, we have
$Te^\al=(\al t^{-1})e^\al = \al_{(-1)}e^\al$.
Using \eqref{twtaa}, \eqref{tvq} and \eqref{twnpr}, we find
\begin{align*}
\d_z Y^M(e^\al,z) &= Y^M(T e^\al,z) = Y^M(\al_{(-1)}e^\al,z) 
\\
&= \nop{  Y^M(\al t^{-1},z)Y^M(e^\al,z) }
- \sum_{j=0}^{N-1} \frac jN (\pi_j\al|\al) z^{-1} Y^M(e^\al,z).
\end{align*}
Therefore
\begin{equation*}
\d_z U^M_\al(z) = -\sum_{j=0}^{N-1} \frac jN (\pi_j\al|\al) z^{-1} U^M_\al(z).
\end{equation*}
Using that $(\pi_j\al|\al) = (\pi_{N-j}\al|\al)$ for $1\le j\le N-1$,
it is easy to see that
\begin{equation}\lbb{xalal}
\sum_{j=0}^{N-1} \frac jN (\pi_j\al|\al) = -b_\al.
\end{equation}
This proves that $U^M_\al = z^{ -b_\al} U^M_\al(z)$
is independent of $z$.

Finally, to prove \eqref{gual}, we apply \eqref{ginv} for $a=e^\al$,
using that $\si(e^\al)=\eta(\al)^{-1}e^{\si\al}$ and that by \eqref{twc5}
we have $E^M_{\si\al}(z) = e^{-2\pi\i\al^M_{(0)}} E^M_\al(e^{2\pi\i}z)$.
\end{proof}

\begin{lemma}\lbb{lyalbe1}
We have
\begin{equation}\lbb{yalbe1}
Y^M(e^\al,z)Y^M(e^\be,w) 
= i_{z,w} f_{\al,\be}(z,w) \, z^{b_\al} w^{b_\be}
\, U^M_\al U^M_\be \, E^M_{\al,\be}(z,w),
\end{equation}
where
\begin{align}
\lbb{falbe}
f_{\al,\be}(z,w) &= 
\prod_{k=0}^{N-1} \bigl(z^{1/N} - \eps^k w^{1/N}\bigr)^{ (\si^k\al|\be) },
\\
\lbb{ealbe}
E^M_{\al,\be}(z,w) &=
\nop{\exp\bigl( \tint Y^M(\al t^{-1},z) +\tint Y^M(\be t^{-1},w) \bigr)}.
\end{align}
\end{lemma}
\begin{proof}
Standard exercise, using the fact that
$
e^Ae^Be^{-A} = e^{\ad A}e^B = e^{[A,B]}e^B
$
for any two operators $A,B$ commuting with $[A,B]$.
\end{proof}

\begin{lemma}\lbb{lyalbe2}
We have
\begin{align}\lbb{yalbe2}
&Y^M(Y(e^\al,z)e^\be,w) 
\\
\notag
&= \ep(\al,\be) B_{\al,\be}^{-1} \, 
i_{w,z} f_{\al,\be}(z+w,w) \, (z+w)^{b_\al} w^{b_\be}
\, U^M_{\al+\be} \, E^M_{\al,\be}(z+w,w),
\end{align}
where
\begin{equation}\lbb{balbe}
B_{\al,\be} = \frac{ f_{\al,\be}(z,w) }{ (z-w)^{(\al|\be)} } 
\Big|_{ z^{1/N}=w^{1/N}=1 }
= N^{ -(\al|\be) } \prod_{k=1}^{N-1} \bigl(1 - \eps^k \bigr)^{ (\si^k\al|\be) }.
\end{equation}
\end{lemma}
\begin{proof}
We use the same argument as in the proof of \leref{lymea}.
First, using \eqref{twcomm}, \eqref{twc4} and \eqref{hteal},
we compute the commutator ($h\in\h_j$, $n\in\frac jN +\Zset$):
\begin{align*}
[h^M_{(n)}, \, & Y^M(Y(e^\al,z)e^\be,w)]
\\
&= \sum_{m=0}^\infty \binom nm w^{n-m} \,
Y^M\bigl( h_{(m)}(Y(e^\al,z)e^\be) ,w \bigr)
\\
&= \sum_{m=0}^\infty \binom nm w^{n-m} 
\bigl( (h|\al)z^m + \de_{m,0} (h|\be) \bigr)
Y^M(Y(e^\al,z)e^\be,w)
\\
&= \bigl( i_{w,z}(z+w)^n (h|\al) + w^n (h|\be) \bigr)
Y^M(Y(e^\al,z)e^\be,w).
\end{align*}
It follows that 
\begin{equation*}
Y^M(Y(e^\al,z)e^\be,w) = U^M_{\al,\be}(z,w) E^M_{\al,\be}(z+w,w),
\end{equation*}
where the operator $U^M_{\al,\be}(z,w)$ satisfies
\begin{equation*}
[h^M_{(n)},U^M_{\al,\be}(z,w)] = 
\de_{n,0}(h_0|\al+\be) \, U^M_{\al,\be}(z,w).
\end{equation*}

Next, we note that (by \eqref{tvq}, \eqref{twtaa}, \eqref{t})
\begin{align*}
Y^M\bigl(Y(e^\al,z)(\be_{(-1)}e^\be),w\bigr)
&= Y^M\bigl(Y(e^\al,z)(Te^\be),w\bigr)
\\
&= (\d_w-\d_z) Y^M(Y(e^\al,z)e^\be,w).
\end{align*}
A similar computation as above, using \eqref{twnpr} for $n=-1$,
shows that for $h\in\h_j$ one has
\begin{multline*}
Y^M\bigl(Y(e^\al,z)(h_{(-1)}e^\be),w\bigr)
= \nop{ Y^M(ht^{-1},w) Y^M(Y(e^\al,z)e^\be,w) }
\\
- \Bigl( z^{-1} i_{w,z} \Bigl(1+\frac zw\Bigr)^{j/N} (h|\al) 
         + \frac jN w^{-1} (h|\be) \Bigr) \,
Y^M(Y(e^\al,z)e^\be,w).
\end{multline*}
Therefore
\begin{align*}
(\d_w-\d_z) & U^M_{\al,\be}(z,w) 
\\
&= -\sum_{j=0}^{N-1} 
\Bigl( z^{-1} i_{w,z} \Bigl(1+\frac zw\Bigr)^{j/N} (\pi_j\be|\al) 
         + \frac jN w^{-1} (\pi_j\be|\be) \Bigr) \, U^M_{\al,\be}(z,w).
\end{align*}
With some more computation, we see that
\begin{align*}
-  \sum_{j=0}^{N-1} z^{-1} i_{w,z} \Bigl(1+\frac zw\Bigr)^{j/N} (\pi_j\be|\al) 
&= i_{w,z} i_{w,z+w} \sum_{k=0}^{N-1}
\frac{ (1/N)\eps^k w^{1/N-1} (\si^k\al|\be) }{ (z+w)^{1/N} - \eps^k w^{1/N} }
\\
&= i_{w,z} \frac{ (\d_w-\d_z) f_{\al,\be}(z+w,w) }{ f_{\al,\be}(z+w,w) }
\end{align*}
On the other hand, by \eqref{xalal}, we have
\begin{equation*}
-\sum_{j=0}^{N-1} \frac jN w^{-1} (\pi_j\be|\be) 
= (\d_w-\d_z) w^{b_\be} / w^{b_\be} .
\end{equation*}
It follows that
\begin{equation*}
U^M_{\al,\be}(z,w) / i_{w,z} f_{\al,\be}(z+w,w) w^{b_\be} 
\end{equation*}
depends only on $z+w$.

Finally, note that 
\begin{align*}
Y^M(Y(e^\al,z)e^\be,w) &= z^{(\al|\be)} \ep(\al,\be) Y^M(e^{\al+\be},w)
+ \text{higher powers of $z$},
\\
\intertext{while} 
i_{w,z} f_{\al,\be}(z+w,w) 
&= z^{(\al|\be)} B_{\al,\be} \, w^{ (\al_0-\al|\be) }
+ \text{higher powers of $z$}.
\end{align*}
Since $(\al_0-\al|\be) = b_{\al+\be} - b_\al - b_\be$,
this completes the proof.
\end{proof}

\begin{corollary}\lbb{cualube}
In any $\si$-twisted $V_Q$-module $M$, one has
\begin{equation}\lbb{ualube}
U^M_\al U^M_\be = \ep(\al,\be) B_{\al,\be}^{-1} \, U^M_{\al+\be},
\qquad \al,\be\in Q.
\end{equation}
\end{corollary}
\begin{proof}
Follows immediately from \eqref{yalbe1}, \eqref{yalbe2} and 
the associativity \eqref{as2}.
\end{proof}

\begin{remark}\lbb{rpfcom}
In the proofs of Lemmas \ref{lyalbe1} and \ref{lyalbe2},
we used only the commutator formulas \eqref{twc3}, \eqref{twc4},
the translation invariance \eqref{twtaa}, and formula \eqref{twnpr}
for $n=-1$, $a=ht^{-1}$, $b=e^\be$.
\end{remark}

\subsection{The Heisenberg Pair $(\hat\h_\si,\G_\si)$}\lbb{shchg}
The results of the previous subsection motivate the following definitions.

The {\em $\si$-twisted current algebra} $\hat\h_\si$ 
consists of all $\si$-invariant elements from 
$\Cset K\oplus\h[t^{1/N},t^{-1/N}]$, where $\si$ acts as
\begin{equation}\lbb{ght}
\si(h t^m) = \si(h) e^{2\pi\i m}t^m, \;\; \si(K)=K,
\qquad h\in\h, \; m\in\tfrac1N\Zset.
\end{equation}
In other words, $\hat\h_\si$ is spanned over $\Cset$ by $K$ and
the elements $ht^m$ such that $h\in\h_j$, $m\in\frac jN +\Zset$.
This is a Lie algebra with bracket \eqref{htht}.

Let $\G = \Cset^\times \times \exp\h_0 \times Q$
be the set consisting of elements $c \, e^h U_\al$ 
($c\in\Cset^\times$, $h\in\h_0$, $\al\in Q$).
We define a multiplication in $\G$ by the formulas:
\begin{align}
\lbb{tig1}
e^h e^{h'} &= e^{h+h'},
\\
\lbb{tig2}
e^h U_\al e^{-h} &= e^{(h|\al)} U_\al,
\\
\lbb{tig3}
U_\al U_\be &= \ep(\al,\be) B_{\al,\be}^{-1} \, U_{\al+\be}.
\end{align}
Then $\G$ is a group.
{}From \eqref{tig3} we get the commutator
\begin{equation}\lbb{tig4}
C_{\al,\be} := 
U_\al U_\be U_\al^{-1} U_\be^{-1} 
= (-1)^{|\al|^2|\be|^2} 
\prod_{k=0}^{N-1} \bigl(-\eps^k \bigr)^{ -(\si^k\al|\be) }.
\end{equation}
We give another expression for $C_{\al,\be}$ which will 
be useful in the sequel:
\begin{equation}\lbb{calbe}
\begin{split}
C_{\al,\be} &= (-1)^{|\al|^2|\be|^2} e^{\pi\i(\al_0|\be)} e^{2\pi\i(\al_*|\be)}
\\
&\qquad\qquad\text{for}\;\; 
\al=\al_0+(1-\si)\al_*, \; \al_0\in\h_0, \al_*\in\h_0^\perp.
\end{split}
\end{equation}

Denote by $Q_\ev$ the sublattice of $Q$ consisting of all even elements, 
i.e., $\al$ such that $|\al|^2 \in 2\Zset$.
\begin{lemma}\lbb{lcentg}
The center $Z(\G)$ of $\G$ consists of all elements of the form
$c \, e^{2\pi\i\la_0} U_{(1-\si)\la}$, where $c\in\Cset^\times$
and $\la\in (Q_\ev)^*$ is such that $\al:=(1-\si)\la\in Q$ and
$\la\in Q^*$ if $\al\in Q_\ev$.
\end{lemma}
\begin{proof}
If $e^{2\pi\i h} U_\al$ is in the center, then \eqref{tig2} implies $\al_0=0$. 
Then $\al=(1-\si)\al_*$ for some uniquely defined $\al_*\in\h_0^\perp$.
Letting $\la=h+\al_*$, we get
$e^{2\pi\i h} U_\al = e^{2\pi\i\la_0} U_{(1-\si)\la}$.
Using \eqref{tig2} and \eqref{calbe}, we see that 
$e^{2\pi\i\la_0} U_{(1-\si)\la}$ commutes with $U_\be$ iff
\begin{equation}\lbb{labe}
(\la|\be)+|\al|^2|\be|^2/2 \in\Zset \quad \text{for $\be\in Q$}. 
\end{equation}
Since
\begin{equation*}
|\al|^2 = ((1-\si)\al_* | (1-\si)\al_*) = 2 (\al_*|\al_*) - 2 (\al_*|\si\al_*)
= 2 (\al_*|\al) = 2 (\la|\al), 
\end{equation*}
equation \eqref{labe} is equivalent
to $\la\in Q^*$ if $\al\in Q_\ev$ and to 
$\la\in (Q_\ev)^*$ if $\al\in Q\setminus Q_\ev$.
\end{proof}

In particular, by \leref{lcentg},
all elements of the form $e^{2\pi\i\al_0} U_{(1-\si)\al}$
($\al\in Q$) are central in $\G$. We let $\G_\si$ be the factor
of $\G$ over the central subgroup 
\begin{equation}\lbb{tig5}
N_\si := \{ \eta(\al) U_{\si\al}^{-1} U_\al e^{ 2\pi\i(b_\al+\al_0) }
\st \al\in Q \} \, .
\end{equation}
Note that $N_\si\cap\Cset^\times = \{1\}$.

We endow $Q$ with the discrete topology so that $\G$ and $\G_\si$ are Lie
groups with a Lie algebra $\Cset\oplus\h_0$. 
By \eqref{tig5}, \eqref{dal} and \eqref{etah0}, we have
\begin{equation}\lbb{tig6}
e^{ 2\pi\i\al } = 1 \quad\text{in $\G_\si$ for $\al\in Q\cap\h_0$}.
\end{equation}
It is easy to see that the connected component of the unit in $\G_\si$
is equal to $\Cset^\times$ times the torus 
\begin{equation}\lbb{tsi}
T_\si := \exp 2\pi\i(\h_0 / Q\cap\h_0).
\end{equation}

The group $\G$ acts on $\hat\h_\si$ by conjugation:
\begin{equation}\lbb{tig7}
(c e^h U_\al) (h't^m + c'K) (c e^h U_\al)^{-1} 
= h't^m + \de_{m,0}(h'_0|\al)K + c'K.
\end{equation}
This action is compatible with the adjoint action of $\Cset\oplus\h_0$
on $\hat\h_\si$ (which is trivial), hence, $(\hat\h_\si,\G)$ is a 
{\em Heisenberg pair\/} in the sense of \cite{FK}. The same is true for 
$(\hat\h_\si,\G_\si)$ because $N_\si$ acts trivially on $\hat\h_\si$.
A module $M$ over $\hat\h_\si$ or over $(\hat\h_\si,\G_\si)$ will be called
{\em restricted\/} if the action of $\h_0$ is diagonalizable
and for any $v\in M$, $(ht^m)v=0$ for $h\in\h$ and sufficiently large 
$m\in\tfrac1N\Zset$.

Now we can summarize the results of \seref{sstwvo} as follows.

\begin{proposition}\lbb{ptvqm}
Any $\si$-twisted $V_Q$-module $M$ is naturally a restricted module over the
Heisenberg pair $(\hat\h_\si,\G_\si)$ of level $1$ 
$($i.e., both $K\in\hat\h_\si$
and $1\in\G_\si$ act as $1)$. Conversely, any 
restricted $(\hat\h_\si,\G_\si)$-module of level $1$ 
can be endowed with the structure of a $\si$-twisted $V_Q$-module.
This establishes an equivalence of the corresponding abelian categories.
\end{proposition}
\begin{proof}
1. Let $M$ be a $\si$-twisted $V_Q$-module. By definition, the
action of $L_0^M$ is diagonalizable with finite-dimensional eigenspaces.
Since $\h_0$ commutes with $L_0^M$, its action is diagonalizable too.
By \eqref{twc3}, the modes $h^M_{(m)}$ ($h\in\h$, $m\in\frac1N\Zset$)
provide a restricted representation of $\hat\h_\si$ of level $1$.
The action of $U_\al\in\G_\si$ is given by the operator $U^M_\al$,
see \eqref{hual}, \eqref{gual}, \eqref{ualube}.

2. Conversely, let $M$ be a restricted $(\hat\h_\si,\G_\si)$-module
of level $1$. Denote the image of $ht^m$ in $\End M$ by $h^M_{(m)}$,
and that of $U_\al$ by $U^M_\al$. This allows us to define the
fields $Y^M(ht^{-1},z)$ and then 
$Y^M(e^\al,z)$ by \eqref{ymeal} ($h\in\h,\al\in Q$),
They satisfy \eqref{ginv}, \eqref{twcomm}; in particular, they are local.
The map $Y^M$ can be extended uniquely to the whole $V_Q$ by applying 
\eqref{as3} repeatedly for $a\in\h t^{-1}$. Note that \eqref{as3}
implies the translation invariance \eqref{twtaa} for $a\in\h t^{-1}$.
Then the proof of \leref{lymea} shows that \eqref{twtaa} holds for
$a=e^\al$, and hence for any $a\in V_Q$.
It follows from \reref{rpfcom} that Lemmas \ref{lyalbe1} 
and \ref{lyalbe2} hold. 
This implies the associativity \eqref{as2}.
By \prref{pdefas}, $M$ is a $\si$-twisted $V_Q$-module.
\end{proof}

\subsection{The Groups $\G^\perp$ and $\G_\si^\perp$}\lbb{sgsgsp}
Before we proceed to the classification of all 
restricted $(\hat\h_\si,\G_\si)$-modules of level $1$,
we need to study the groups $\G$ and $\G_\si$ in more detail.

Let $\G^\perp \subset \G$ be the subgroup of $\G$ consisting
of all $c \, U_\al$ with $c\in\Cset^\times$, $\al\in Q\cap\h_0^\perp$. 
Clearly, the centralizer of $\hat\h_\si$ in $\G$ equals
$\exp\h_0 \times \G^\perp$ (cf.\ \eqref{tig7}). 
In other words, $\G^\perp$ is the outer centralizer 
of the torus $\exp\h_0$ in~$\G$. 

Denote by $\G_\si^\perp$ 
the image of $\G^\perp$ in $\G_\si$.
It can be described as the factor of $\G^\perp$ over the
central subgroup (cf.\ \eqref{tig5}, \eqref{dal})
\begin{equation}\lbb{tigp}
N_\si^\perp := N_\si\cap\G^\perp
= \{ \eta(\al)(-1)^{|\al|^2} U_{\si\al}^{-1} U_\al
\st \al\in Q\cap\h_0^\perp \} \, . 
\end{equation}
The centralizer of $\hat\h_\si$ (and of $T_\si$) in $\G_\si$ is equal to
$T_\si \times \G_\si^\perp$.
Notice that $\G_\si^\perp$ is a central extension (by $\Cset^\times$)
of the finite abelian group $(Q\cap\h_0^\perp) / (1-\si)(Q\cap\h_0^\perp)$.

\begin{definition}\lbb{dzsi}
\textup{(i)}
Let $P_\si$ be the set of all $\la$ that appear in \leref{lcentg}, i.e.,
the set of all $\la\in (Q_\ev)^*$ such that $(1-\si)\la\in Q$ 
and $\la\in Q^*$ if $(1-\si)\la\in Q_\ev$.
Note that $P_\si$ is a sublattice of $(Q_\ev)^*$
containing $Q$.

\textup{(ii)}
Let $Q_\si = (1-\si)P_\si \subset Q$.

\textup{(iii)}
Let $Z_\si = P_\si/Q$ be the subgroup of $((Q_\ev)^*/Q)^\si$ consisting of
classes $\la+Q$ such that $\la\in Q^*$ if $(1-\si)\la\in Q_\ev$.
In particular, when the lattice $Q$ is even, $Z_\si = (Q^*/Q)^\si$ 
is the group of $\si$-invariant elements in $Q^*/Q$.
\end{definition}

Similarly to \leref{lcentg}, we can describe the centers of $\G_\si$, 
$\G^\perp$ and~$\G_\si^\perp$.

\begin{lemma}\lbb{lcgg}
\textup{(i)}
$Z(\G_\si) \simeq Z(\G)/N_\si \simeq \Cset^\times \times Z_\si$.

\textup{(ii)}
$Z(\G^\perp) 
= \{ c \, U_\al \st c\in\Cset^\times, \al\in Q_\si \}
\simeq \Cset^\times \times Q_\si$.

\textup{(iii)}
$Z(\G_\si^\perp) \simeq Z(\G^\perp) / N_\si^\perp \simeq
\Cset^\times \times Q_\si / (1-\si)(Q\cap\h_0^\perp)$.
%
%
\end{lemma}
\begin{proof}
(i) follows from \leref{lcentg}, \eqref{calbe} and the fact that
$N_\si\cap\Cset^\times = \{1\}$.

(ii) A similar argument as in the proof of \leref{lcentg} shows that
the center of $\G^\perp$ consists of all elements of the form
$c \, U_{(1-\si)\la}$, where $c\in\Cset^\times$
and $\la\in (Q_\ev\cap\h_0^\perp)^*$ is such that $\al:=(1-\si)\la\in Q$ and
$\la\in (Q\cap\h_0^\perp)^*$ if $\al\in Q_\ev$.
Next, we use the following lemma.

\begin{lemma}\lbb{lproj*}
Let $\h'$ be a subspace of $\h$ such that the restriction of the
bilinear form on it is nondegenerate. Denote by $\pi'$ the orthogonal
projection of $\h$ onto $\h'$. Then for any lattice $L\subset\h$, one has
$L^*\cap\h' = (\pi' L)^{*'}$, where $*'$ means that the dual is taken
in $\h'$.
\end{lemma}
\begin{proof}
Follows from the fact that $(h|\al)=(h|\pi'\al)$
for $h\in\h'$, $\al\in L$.
\end{proof}

Now, by \leref{lproj*}, $(Q\cap\h_0^\perp)^* = \pi_\perp (Q^*)$, 
and similarly for $Q_\ev$, 
where $\pi_\perp$ is the orthogonal projection from $\h$ to $\h_0^\perp$.
Noting that $(1-\si) \pi_\perp = 1-\si$ completes the proof of part (ii).

Part (iii) follows from part (ii), \eqref{calbe} and the fact that
$N_\si^\perp\cap\Cset^\times = \{1\}$.
\end{proof}

\begin{corollary}\lbb{ccgg}
There is a natural exact sequence
\begin{equation}\lbb{zseq}
1 \to \exp 2\pi\i (Q^*\cap\h_0) \to Z(\G) \overset{p}\to Z(\G^\perp) \to 1 
\, ,
\end{equation}
where the homomorphism $p$ is given by
\begin{equation}\lbb{zseq2}
c \, e^{2\pi\i\la_0} U_{(1-\si)\la} \overset{p}\mapsto c \, U_{(1-\si)\la}
\, , \qquad c\in\Cset^\times \, , \; \la\in P_\si \, .
\end{equation}
The sequence \eqref{zseq} splits, so we have a non-canonical isomorphism
$Z(\G) \simeq \exp 2\pi\i (Q^*\cap\h_0) \times Z(\G^\perp)$.
%
\end{corollary}
\begin{proof}
Clearly, the kernel of $p$ consists of $e^{2\pi\i\la_0}$ with $\la\in Q^*$ 
such that $(1-\si)\la = 0$.
Let $s\colon \pi_\perp(P_\si) \to P_\si$ 
be a linear section of the projection $\pi_\perp$.
For $\la'\in\pi_\perp(P_\si)$, let $\la=s(\la')\in P_\si$.
Then $\la'=\pi_\perp(\la)$ and $(1-\si)\la = (1-\si)\la'$.
The map 
$c \, U_{(1-\si)\la'} \mapsto c \, e^{2\pi\i \la_0} U_{(1-\si)\la}$
is a splitting of \eqref{zseq}.
\end{proof}

Since $\G_\si^\perp$ is a central extension of a finite abelian
group, its representations are completely reducible and
the irreducible ones are classified by the characters of $Z(\G_\si^\perp)$.
We will consider only representations on which $1\in\G_\si^\perp$ acts
as the identity operator. The irreducible ones are classified by the
finite abelian group $Q_\si / (1-\si)(Q\cap\h_0^\perp)$. 
All of them have the same
dimension $d(\si)$, which satisfies
\begin{equation}\lbb{ddef1}
d(\si)^2 = |\G_\si^\perp/ Z(\G_\si^\perp)| = |(Q\cap\h_0^\perp) / Q_\si| .
\end{equation}
\begin{definition}\lbb{ddefect}
The non-negative integer 
$d(\si)$ 
is called the {\em defect\/} of $\si$ (cf.~\cite{KP}).
\end{definition}

\subsection{Representations of\/ $(\hat\h_\si,\G_\si)$}\lbb{srephg}
In this subsection we show that the category of all restricted
$(\hat\h_\si,\G_\si)$-modules of level $1$ is semisimple, and we classify
the irreducible ones.

Let $\hat\h_\si^-$ (resp.\ $\hat\h_\si^+$) be the subalgebra of $\hat\h_\si$
consisting of all elements $ht^m$ with $m>0$ (resp.\ $m<0$).
It is well known (see, e.g., \cite{K1}) that any restricted 
$\hat\h_\si$-module $M$ of level $1$ is induced from its vacuum
subspace
\begin{equation}\lbb{omm}
\Om_M = \{ v\in M \st \hat\h_\si^- v=0 \}.
\end{equation}
More precisely, 
\begin{equation}\lbb{omm2}
M \simeq \Ind^{\hat\h_\si}_{\h_0\oplus\hat\h_\si^-} \Om_M
\simeq S(\hat\h_\si^+)\tt\Om_M .
\end{equation}
The subalgebra $\h_0$ acts on $\Om_M$ diagonally, and the 
$\hat\h_\si$-module $M$ 
is completely reducible (it is irreducible iff $\dim\Om_M=1$).

Now assume that $M$ is a restricted $(\hat\h_\si,\G_\si)$-module of level $1$.
It follows from \eqref{tig7} that $\Om_M$ is a $\G_\si$-module.
By definition, $1\in\G_\si$ acts as $1$ and the torus $T_\si$
acts diagonally (cf.\ \eqref{tsi}). We will call such $\G_\si$-modules
{\em restricted\/}.

If $\Om$ is a restricted $\G_\si$-module, 
it has a compatible $\h_0$-action, because $\h_0$
is the Lie algebra of the torus $T_\si$.
We let $\hat\h_\si^-$ act trivially on $\Om$
and form the induced $\hat\h_\si$-module 
$M(\Om) = \Ind^{\hat\h_\si}_{\h_0\oplus\hat\h_\si^-}\Om$.
Using \eqref{tig7}, we can extend the action of $\G_\si$
from $\Om$ to $M(\Om)$. Then $M(\Om)$ becomes a 
restricted $(\hat\h_\si,\G_\si)$-module of level $1$.

\begin{proposition}\lbb{pmomm}
The functors $M\mapsto\Om_M$ and 
$\Om\mapsto M(\Om)$
establish an equivalence of abelian categories between 
the category of restricted $(\hat\h_\si,\G_\si)$-modules 
of level $1$ and the category of 
restricted $\G_\si$-modules.
\end{proposition}

Therefore we are left with describing restricted $\G_\si$-modules.

\begin{proposition}\lbb{padmgg}
Any restricted $\G_\si$-module is completely reducible and is determined
by the action of the center of $\G_\si$. 
Isomorphism classes of restricted irreducible $\G_\si$-modules are
parameterized by the {\rm(}finite{\rm)} set $Z_\si$.
\end{proposition}

Let $\Om$ be a restricted $\G_\si$-module. 
For $\mu\in\h_0$, we denote by $\Om_\mu$ the weight $\mu$
subspace of $\Om:$
\begin{equation}\lbb{ommu}
\Om_\mu := \{ v\in\Om \st e^h v = e^{(h|\mu)} v 
\;\;\text{for}\;\; h\in\h_0 \}. 
\end{equation}
Then
$(\mu|\al)\in\Zset$ for $\al\in Q\cap\h_0$, i.e., 
$\mu\in(Q\cap\h_0)^* = \pi_0(Q^*)$ by \leref{lproj*}.

\begin{lemma}\lbb{ladmgg1}
\textup{(i)}
$U_\al\Om_\mu = \Om_{\mu+\pi_0\al}$
for all $\al\in Q$, $\mu\in\pi_0(Q^*)$.
In particular,
the subgroup $\G_\si^\perp \subset \G_\si$ 
preserves each $\Om_\mu$.

\textup{(ii)}
Let 
$\Om_{\mu+\pi_0(Q)} = \sum_{\al\in Q} \Om_{\mu+\pi_0\al} \subset \Om$.
Then $\Om_{\mu+\pi_0(Q)}$ is a $\G_\si$-submodule of~$\,\Om$.

\textup{(iii)}
The $\G_\si$-module $\Om_{\mu+\pi_0(Q)}$ is irreducible 
if and only if the $\G_\si^\perp$-module
$\Om_\mu$ is irreducible. 
\end{lemma}
\begin{proof}
(i) 
It follows from \eqref{tig2} that $U_\al\Om_\mu \subset \Om_{\mu+\pi_0\al}$
for any $\al\in Q$. Since by \eqref{tig3} $U_\al^{-1}$ is proportional
to $U_{-\al}$\,, we get $U_\al\Om_\mu = \Om_{\mu+\pi_0\al}$\,.

(ii) follows from (i) and the definition of $\G_\si$ 
(see \eqref{tig1}--\eqref{tig3}).

(iii)
First, let $\Om_\mu$ be an irreducible $\G_\si^\perp$-module.
Assume that $\La$ is a nontrivial $\G_\si$-submodule of $\Om_{\mu+\pi_0(Q)}$.
Using the action of $T_\si$, we can write 
$\La = \sum_{\al\in Q} \La_{\mu+\pi_0\al}$
where each $\La_{\mu+\pi_0\al} \subset \Om_{\mu+\pi_0\al}$\,.
Moreover, $\La_{\mu+\pi_0\al} = U_\al \La_\mu$.
In particular, $\La_\mu$ is a $\G_\si^\perp$-submodule of $\Om_\mu$.
But $\Om_\mu$ is irreducible; hence,
$\La_\mu = \Om_\mu$ and $\La = \Om_{\mu+\pi_0(Q)}$.

Conversely, assume that the $\G_\si^\perp$-module $\Om_\mu$ 
is not irreducible. Since $\G_\si^\perp$ is a central extension 
of a finite abelian group, its representations are completely reducible.
If $\Om_\mu = \bigoplus_i L_i$ as a $\G_\si^\perp$-module, let
$L^i  = \sum_{\al\in Q} U_\al L_i$.
Then each $L^i$ is a $\G_\si$-submodule of $\Om_{\mu+\pi_0(Q)}$
and $\Om_{\mu+\pi_0(Q)} = \bigoplus_i L^i$ as a $\G_\si$-module.
\end{proof}
\begin{remark}\lbb{romind}
$\Om_{\mu+\pi_0(Q)}$ is isomorphic to the induced module 
$\Ind^{\G_\si}_{T_\si \times \G_\si^\perp} \Om_\mu$.
\end{remark}

{}From \leref{ladmgg1}(iii) and its proof, we see that the
$\G_\si$-module $\Om_{\mu+\pi_0(Q)}$ is completely reducible.
Since
$\Om = \bigoplus_{ [\mu]\in\pi_0 (Q^*)/\pi_0(Q) } \Om_{[\mu]}$,
it is also completely reducible.

Now let $\Om$ be an irreducible $\G_\si$-module.
Then $\Om=\Om_{\mu+\pi_0(Q)}$ for some $\mu\in\pi_0 (Q^*)$
and the $\G_\si^\perp$-module $\Om_\mu$ is irreducible.
Any irreducible $\G_\si^\perp$-module is completely determined
by the action of the center $Z(\G_\si^\perp)$.
Let $\ze\colon Z(\G_\si^\perp) \to \Cset^\times$
be the central character of $\Om_\mu$.

We can view $\Om$ as a $\G$-module on which $N_\si$ acts trivially,
and similarly $\Om_\mu$ as a $\G^\perp$-module with a trivial action
of $N_\si^\perp$. Recall that, by \leref{lcgg}(iii), 
$Z(\G_\si^\perp) \simeq Z(\G^\perp) / N_\si^\perp$,
so we can extend $\ze$ to a character of $Z(\G^\perp)$.

If $\mu'=\mu+\pi_0\al$ for some $\al\in Q$, then
$\Om_{\mu'}=U_\al\Om_\mu$ and $\Om_{\mu+\pi_0(Q)}=\Om_{\mu'+\pi_0(Q)}$.
For $v\in\Om_\mu$, $U_\be\in Z(\G^\perp)$, $\be\in Q_\si$
(see \leref{lcgg}(ii)), we have:
$U_\be U_\al v = C_{\al,\be}^{-1} U_\al U_\be v 
= C_{\al,\be}^{-1} \ze(U_\be) U_\al v$
where $C_{\al,\be}$ is given by \eqref{tig4}.
Hence, two pairs $(\mu,\ze)$ and $(\mu',\ze')$
correspond to the same irreducible $\G_\si$-module if and only if
they are related by:
\begin{equation}\lbb{muze}
\mu'=\mu+\pi_0\al \,, \quad 
\ze'(U_\be) = C_{\al,\be}^{-1} \ze(U_\be) \,, \qquad
\al\in Q \,, \; \be\in Q_\si \,.
\end{equation}

For $\la\in P_\si$ the element 
$e^{2\pi\i\la_0} U_{(1-\si)\la} \in Z(\G)$
acts on $\Om_\mu$ as the scalar 
\linebreak
$e^{2\pi\i(\la_0|\mu)} \ze(U_{(1-\si)\la})$
(cf.\ Lemmas \ref{lcentg} and \ref{lcgg}(ii)).
Using \coref{ccgg} and \leref{lproj*},
it is easy to see that the action of $Z(\G)$ on $\Om$
determines uniquely the equivalence class of $(\mu,\ze)$ under
\eqref{muze}, and hence it determines the isomorphism
class of the $\G_\si$-module $\Om$.
Conversely, different pairs $(\mu,\ze)$ give rise to different
actions of $Z(\G)$ on the corresponding modules $\Om_{\mu+\pi_0(Q)}$.

This completes the proof of \prref{padmgg}.

\subsection{Classification of $\si$-Twisted $V_Q$-Modules}\lbb{sclastwm}
Combining Propositions \ref{ptvqm}, \ref{pmomm} and \ref{padmgg},
we obtain the main result of the paper.
\begin{theorem}\lbb{ttwvqm}
The category of $\si$-twisted $V_Q$-modules is a semisimple abelian category
with finitely many isomorphism classes of simple objects,
parameterized by the set $Z_\si$. 
\end{theorem}

\begin{remark}\lbb{racth}
The irreducible $\si$-twisted $V_Q$-module
corresponding to $\la+Q\in Z_\si$
is isomorphic as an $\hat\h_\si$-module to 
$S(\hat\h_\si^+) \tt e^{\la_0}\Cset[\pi_0 Q] \tt\Cset^{d(\si)}$,
where $\Cset$ carries the zero action and
$d(\si)$ is the defect of $\si$.
\end{remark}


\end{document}